\documentclass[12pt]{article} \usepackage{amssymb} \oddsidemargin 0pt
\newtheorem{theo}{Theorem}[section]
\newtheorem{prop}[theo]{Proposition}
\newtheorem{lemma}[theo]{Lemma}
\newtheorem{coro}[theo]{Corollary}
\newtheorem{conj}[theo]{Conjecture}
\newtheorem{prob}[theo]{Problem}
\newcommand{\ignore}[1]{}
\def\square{\vrule height6pt width7pt depth1pt}
\def\endpf{\hfill\square\bigskip}
\begin{document}

\title{A Tur\'{a}n Type Problem Concerning the Powers of the Degrees
of a Graph}
\author{Yair Caro \thanks{e-mail: yairc@macam98.ac.il}\\
and\\ Raphael Yuster \thanks{e-mail: raphy@research.haifa.ac.il}\\ Department of
Mathematics\\ University of Haifa-ORANIM, Tivon 36006, Israel.}
\date{}
\maketitle

\begin{abstract}
For a graph $G$ whose degree sequence is $d_{1},\ldots ,d_{n}$, and for a
positive integer $p$, let $e_{p}(G)=\sum_{i=1}^{n}d_{i}^{p}$. For a fixed
graph $H$, let $t_{p}(n,H)$ denote the maximum value of $e_{p}(G)$ taken over
all graphs with $n$ vertices that do not contain $H$ as a subgraph. Clearly,
$t_{1}(n,H)$ is twice the Tur\'{a}n number of $H$. In this paper we consider
the case $p>1$. For some graphs $H$ we obtain exact results, for some others
we can obtain asymptotically tight upper and lower bounds, and many
interesting cases remain open.
\end{abstract}

\section{Introduction}

All graphs considered here are finite, undirected, and have no loops or
multiple edges. For the standard graph-theoretic notations the reader is
referred to \cite{Bo}. For a graph $G$ whose degree sequence is $
d_{1},\ldots ,d_{n}$ let $e_{p}(G)=\sum_{i=1}^{n}d_{i}^{p}$. Clearly,
$e_{1}(G)=2e(G)$. Recently, several papers were published concerning the
problem of maximizing $e_{2}(G)$ over all graphs having $n$ vertices and $m$
edges. See, e.g., \cite{By,CaYu,Ol,De,PePeSt}. In this line of research no
restriction is imposed on the structure of $G$. Along the spirit of Tur\'{a}n
Theory we consider the problem of finding the maximum of $e_{p}(G)$ over the
class of graphs which contain no copy of prescribed forbidden subgraphs. For
a fixed graph $H$, let $t_{p}(n,H)$ denote the maximum value of $e_{p}(G) $
taken over all graphs with $n$ vertices that do not contain $H $ as a
subgraph. Clearly, $t_{1}(n,H)=2t(n,H)$ where $t(n,H)$ is the Tur\'{a}n
Number of $H$.

The Tur\'an number $t(n,H)$ is one of the most studied parameters in Graph
Theory. Many interesting and non-trivial results give either exact values or
asymptotically tight upper and lower bounds for $t(n,H)$. For example, the
classic result of Tur\'an (from which Tur\'an Theory has emerged) determines
$t(n,K_p)$ for all $n$ and $p$. There are still many open problems, even when
$H$ is a rather simple graph. For example, when $H$ is a tree with $k$
vertices it is conjectured that $t(n,H)=(k/2-1)n(1+o(1))$, and when $k-1$
divides $n$ the conjecture is that $t(n,H)=(k/2-1)n$. The lower bound is
obtained by taking $n/(k-1)$ vertex-disjoint copies of $K_{k-1}$. The upper
bound would follow if one can prove the famous conjecture of Erd\H{o}s and
S\'os \cite{Er0}, which states that graphs with $(k/2-1)n+1$ edges contain
every tree with $k$ vertices. This conjecture is known to hold if $G$ is
$C_4$-free \cite{SaWo} or if the tree has a vertex adjacent to at least
$(k-2)/2$ leaves \cite{Si}.

In many cases, the extremal graphs with respect to $t(n,H)$ tend to be
regular or almost regular. That is, the $k$'th central moment of the degree
sequence is either zero or very small. If we wish to investigate highly
non-regular $H$-forbidden graphs, then just counting the number of edges does
not suffice. If we wish to maximize the second central moment of the degree
sequence of $H$-forbidden graphs, then the parameter $t_2(n,H)$ is the
correct measure. Likewise, for the $p$'th central moment the parameter
$t_p(n,H)$ is the suitable one. In this paper we consider $t_p(n,H)$ for $p >
1$. For some graphs $H$ we are able to give exact or near-exact results,
while for others the problem remains open.

Our first result shows that in both the quadratic and cubic cases, when $H=K_k$, the extremal graph that yields
$t_p(n,H)$ is exactly the same graph that yields $t_1(n,H)$, namely the
Tur\'an Graph $T(n,k)$.

\begin{theo}
\label{t11} Let $k > 2$ be a positive integer, and let $p=1,2,3$. Then, 
$t_p(n,K_k)=e_p(T(n,k))$, where $T(n,k)$ is the Tur\'an Graph.
\end{theo}
Theorem \ref{t11} is sharp in the sense that for $p \geq 4$, $t_p(n,K_k)$ is {\em not}
obtained by the Tur\'an graph (but is obtained by another non-balanced complete $(k-1)$-partite graph).
Note: In the original version of this paper Theorem \ref{t11} was mistakenly stated for all $p$.
This was also observed by Pikhurko and by Schelp.

Let $P_{k}$ denote the path with $k$ vertices. Faudree and Schelp \cite{FaSc}
characterized the extremal graphs that yield $t(n,P_{k})$. Let $r\equiv
n\bmod(k-1)$. An extremal graph giving $t(n,P_{k})$ is obtained by taking
$\lfloor n/(k-1)\rfloor $ vertex-disjoint copies of cliques of order $k-1$
and, if $r\neq 0$, another clique $K_{r}$ on the remaining vertices. Hence,
$ t(n,P_{k})={{k-1} \choose {2}}\lfloor n/(k-1)\rfloor +{r \choose 2}$.
These graphs are far from optimal when considering $t_{p}(n,P_{k})$ when
$p>1$, since they have small maximum degree.

Our next theorem determines $t_{p}(n,P_{k})$ for $n$ sufficiently large (for
small values of $n$ there are some disturbances). In order to describe this
theorem we define the graph $H(n,k)$ for $n\geq k\geq 4$ as follows. The
vertex set of $H$ is composed of two parts $A$ and $B$ where $|B|=\lfloor
k/2\rfloor -1$ and $|A|=n-|B|$. $B$ induces a complete graph, and $A$ induces
an independent set when $k$ is even, or a single edge plus $|A|-2$ isolated
vertices when $k$ is odd. All possible edges between $A$ and $B$ exist.

\begin{theo}
\label{t12} Let $k\geq 4$, let $p\geq 2$ and let $n>n_{0}(k)$. Then, $H(n,k)$
contains no copy of $P_{k}$ and $t_{p}(n,P_{k})=e_{p}(H(n,k))$. Furthermore,
$H(n,k)$ is the unique extremal graph.
\end{theo}

Note that, trivially, $t_{p}(n,P_{2})=0$ and $t_{p}(n,P_{3})=n$ when $n$ is
even and $t_{p}(n,P_{3})=n-1$ when $n$ is odd (by taking a maximum matching
on $n$ vertices). Also note that when $n$ is small compared to $k$, the graph
$H(n,k)$ is \emph{not} the extremal graph. As an extreme example note that
$t_{p}(k-1,P_{k})=(k-1)(k-2)^{p}$ and is obtained by $K_{k-1}$. A close
examination of the proof of Theorem \ref{t12} shows that the value of
$n_{0}(k)$ in the statement of the theorem is $O(k^{2})$. Another thing to
note is that, as long as $p\geq 2,$ the actual value of $p$ is insignificant.

Let $C^{\ast }$ be the family of even cycles. It is an easy exercise to show
that any graph with more than $\lfloor 3(n-1)/2\rfloor $ edges contains an
even cycle. This bound is sharp and there are exponentially many extremal
graphs \cite{Bo}. In fact, the extremal graphs can be constructed recursively
as follows. For $n=1$ take a single point. For $n=2$ take a single edge. If
$n>2$ we construct graphs with no even cycles and with 
$ \lfloor 3(n-1)/2\rfloor $ edges as follows. Let $G$ be any such extremal
graph with
$n-2$ vertices. Pick any vertex $x$ of $G$ and add to $G$ two new vertices
$a,b$. Now add a triangle on $x,a,b$. The resulting graph has $n$ vertices
$e(G)+3=\lfloor 3(n-1)/2\rfloor $ edges, and no even cycle. Notice that the
\emph{Friendship Graph} $F_{n}$ is one of the extremal graphs. $ F_{n}$ is
defined as follows. Take a star with $n$ vertices and add a maximum matching
on the set of leaves. Thus, $F_{n}$ has exactly $ n-1+\lfloor (n-1)/2\rfloor
=\lfloor 3(n-1)/2\rfloor $ edges, and no even cycle. Note that when $n$ is
odd, $e_{2}(F_{n})=(n-1)^{2}+4(n-1)=n^{2}+2n-3$ and when $n$ is even
$e_{2}(F_{n})=(n-1)^{2}+4(n-2)+1=n^{2}+2n-6$. Our next theorem shows that,
unlike the Tur\'{a}n case, there is only one extremal graph giving
$t_{2}(n,C^{\ast })$, and it is $F_{n}$. (Notice the natural extension of the
definition of $t_{p}$ to families of graphs).

\begin{theo}
\label{t13} For $n$ sufficiently large, $t_2(n,C^*)=e_2(F_n)$ and $F_n$ is
the unique extremal graph.
\end{theo}

We mention that Theorem \ref{t13} also holds for $p > 2$, but the proof is
rather technical and we omit it.

The rest of this paper is organized as follows. In Section 2 we consider
complete graphs and prove Theorem \ref{t11}. In Section 3 we consider paths
and prove Theorem \ref{t12}. Some other acyclic graphs $H$ for which 
$t_p(n,H)$ can be determined are handled in Section 4. In Section 5 we prove
Theorem \ref{t13} and also asymptotically determine $t_k(n,K_{k,k})$. The
final section contains some concluding remarks and open problems.

\section{Complete graphs}

In order to prove Theorem \ref{t11} we need the following theorem of
Erd\H{o}s \cite{Er} that characterizes the maximal degree sequences of
graphs without a $K_k$.

\begin{lemma}[Erd\H{o}s \protect\cite{Er}]
\label{l21} Let $G=(V,E)$ be a graph without a $K_k$. Then, there is a $(k-1)
$-partite graph $G^{\prime}=(V,E^{\prime})$ such that for every $v \in V$, 
$d_G(v) \leq d_{G^{\prime}}(v)$. \endpf
\end{lemma}

If $G$ and $G^{\prime}$ are as in Lemma \ref{l21} then, clearly, $e_p(G) \leq
e_p(G^{\prime})$ for all $p \geq 1$. Thus, the following corollary
immediately follows from Lemma \ref{l21}:

\begin{coro}
\label{c21} For every $n\geq k-1\geq 1$ there exists a complete
$(k-1)$-partite
graph $G$ with $n$ vertices such that $t_{p}(n,K_{k})=e_{p}(G)$.
\endpf
\end{coro}
\textbf{Proof of Theorem \ref{t11}:}\thinspace\ It suffices to show that the
complete $(k-1)$-partite graph $G$ in Corollary \ref{c21} is the Tur\'{a}n
Graph $T(n,k)$, in the cases $p=2,3$. For $k=2$ this is trivial. Assume therefore that $k\geq 3$.
It suffices to show that if $G^{\prime }$ is any complete $(k-1)$-partite
graph that has (at least) two vertex classes $X$ and $Y$ with $|X|-|Y|>1$
then the complete $(k-1)$-partite graph $G^{\prime \prime }$ obtained from
$G^{\prime }$ by transferring a vertex from $X$ to $Y$ has $e_{p}(G^{\prime
\prime })>e_{p}(G^{\prime })$. Indeed, putting $|X|=x$ and $|Y|=y$ we have \[
e_{p}(G^{\prime \prime })-e_{p}(G^{\prime
})=(y+1)(n-y-1)^{p}+(x-1)(n-x+1)^{p}-y(n-y)^{p}-x(n-x)^{p}>0 \] where the
last inequality may be verified using standard (although tedious) calculus,
and the facts that $n\geq x+y$ and $x-y-1>0$. For example, if $p=2$ the
expression in the middle of the last inequality is equivalent to the
expression $(x-y-1)(n+3(n-x-y))$. \endpf

Theorem \ref{t11} is not true for $p \geq 4$. This can already be seen by the fact that the complete bipartite graph
$G=K_{\lfloor n/2-1 \rfloor, \lceil n/2+1 \rceil}$
has $e_4(G) > e_4(T(n,3))$.

Let $K_k^{\prime}$ be the graph obtained from $K_k$ by adding a new vertex of
degree one, connected to one of the original vertices. It is not difficult to
show that $t_p(n,K_k^{\prime})=t_p(n,K_k)$ for $n> n(k)$ (assuming $k \geq 3$
and $p \geq 2$). Indeed, we can state this more generally.

\begin{prop}
\label{p22} Let $H$ be a vertex-transitive graph with at least two edges.
Let $H^{\prime}$ be obtained from $H$ by adding a new vertex of degree one
connected to one of the original vertices. Then, if $p \geq 2$,
$t_p(n,H^{\prime})=t_p(n,H)$ for $n> n(p,H)$.
\end{prop}
\textbf{Proof:}\, Clearly, $t_p(n,H^{\prime}) \geq t_p(n,H)$ since 
$H^{\prime}$ contains $H$. Now assume that equality does not hold. Let $G$ be
an $n$-vertex graph having a copy of $H$ as a subgraph, but having no
$H^{\prime}$ as a subgraph, and having $e_p(G)=t_p(n,H^{\prime}) > t_p(n,H)$.
Since $H$ is vertex-transitive, the set of vertices of every copy of $H$ in
$G$ is disconnected from the other vertices of $G$, since otherwise we
would have an $H^{\prime}$. Thus, if $t$ is the number of copies of $H$ in
$G$ we have $e_p(G)\leq th(h-1)^p+t_p(n-ht,H)$. However,
$t_p(n,H)=\Omega(n^p)$ as can be seen by the star $S_n$ which has no copy of
$H$ (recall that $H$ is vertex-transitive with at least two edges), and
$e_p(S_n)=\Omega(n^p)$, thus for $n$ sufficiently large,
$th(h-1)^p+t_p(n-ht,H)$ is maximized when $t=0$. Consequently, $e_p(G) \leq
t_p(n,H)$, a contradiction. \endpf

\section{Paths}

In order to determine $t_p(n,P_k)$ it is useful to have an upper bound on the
maximum number of edges possible in a graph not containing $P_k$. The
following lemma, which is a theorem of Faudree and Schelp, determines the
Tur\'an number for paths.

\begin{lemma}[Faudree and Schelp \protect\cite{FaSc}]
\label{l31} Let $k > 1$ and let $n > 0$. Let $r=\lfloor n/(k-1) \rfloor$ and
let $s \equiv n \bmod (k-1)$ where $0 \leq s < k-1$. Then $t(n,P_k) = r{{k-1}
\choose 2}+ {s \choose 2}$. \endpf
\end{lemma}

In fact, Faudree and Schelp also characterized the extremal Tur\'an graphs.
The graph composed of $r$ vertex-disjoint cliques of order $k-1$ plus an
additional clique of order $s$ is extremal (sometimes, however, it is not the
only extremal graph). It would be somewhat more convenient to use the
following less accurate upper bound for $t(n,P_k)$, that is always at least
as large as the value in Lemma \ref{l31}

\begin{coro}
\label{c32} If $G$ has $n$ vertices and is $P_k$-free then $e(G) \leq
n(k-2)/2$. \endpf
\end{coro}

We also need a lemma bounding $e_p(G)$ for $n$-vertex graphs $G$ that have
linearly many edges, and have maximum degree $\Theta(n)$.

\begin{lemma}
\label{l33} Let $p \geq 2$ be an integer, let $0.5 < \alpha \leq 1$ and let
$t > \alpha$ be real. Let $G$ be an $n$-vertex graph with $\Delta(G) \leq
\alpha n$ and with at most $tn$ edges. Then: \[ e_p(G) \leq
\frac{t}{\alpha}(\alpha n)^p + o(n^p). \]
\end{lemma}
\textbf{Proof:}\, Consider the degree sequence of $G$, denoted
$\{d_1,\ldots,d_n\}$. It is a sequence of $n$ nonnegative integers whose sum
is at most $2tn$ and whose elements do not exceed $\alpha n$. If we
``forget'' that this sequence is graphic then $d_1^p+\ldots+d_n^p$ is,
obviously, at most \[ \frac{2tn}{\alpha n}(\alpha n)^p =
\frac{2t}{\alpha}(\alpha n)^p. \] However, the sequence is graphic. This
means, for example, that if there is a vertex of degree, say, $\alpha n$,
then there are at least $\alpha n$ nonnegative elements in the sequence. In
fact, for any fixed $\beta \leq \alpha$ there are at most $t/\beta+o(1)$
vertices with degree at least $ \beta n$, and if this happens, then the
other degrees are all at most $ t/\beta+o(1)$ which is constant (and hence
have no significant contribution to $e_p(G)$). By the convexity of the
polynomial $x^p$, the optimal situation is obtained by taking vertices with
degree $\alpha n$ as many as possible (there are at most $t/\alpha+o(1)$ such
vertices), and this forces the other vertices (except maybe one) to have
constant degree. Hence, \[ e_p(G) \leq \frac{t}{\alpha}(\alpha n)^p+o(n^p).
\] \endpf

Before we prove Theorem \ref{t11} we need to dispose of the special case
$k=5 $, since this value causes technical difficulties in the proof.

\begin{lemma}
\label{l340} For $n \geq 12$, $t_p(n,P_5)=e_p(H(n,5))$. Furthermore, $H(n,5)$
is the unique extremal graph.
\end{lemma}
\textbf{Proof:}\, First note that for $r \geq 3$, $H(r,5)$ is a star with $r$
vertices with an additional edge connecting two of its leaves. Assume $G$ is
an $n$-vertex graph having no $P_5$ and $t_p(n,P_5)=e_p(G)$. Trivially, $G$
may only contain cycles of length 3 or 4. In fact, the vertices of every
4-cycle must induce a component of $G$ isomorphic to $K_4$, and for any
3-cycle, the connected component to which it belongs must be an $H(r,5)$ for
some $r \geq 3$. Thus, the components of $G$ are either $K_4$'s or $H(r,5)$'s
(there may be several with distinct values of $r$), or trees. Trivially,
every tree $T$ with $r$ vertices has $e_p(T) \leq e_p(S_r)$ where $S_r$ is
the $r$-vertex star. Similarly, $e_p(S_r) \leq e_p(H(r,5))$ since $H(r,5)$
contains $S_r$ (if $r=1$ or $r=2$ we define $H(1,5)=S_1$ and $H(2,5)=S_2$).
Thus, we may assume that every component is either a $K_4$ or an $H(r,5)$.
Another trivial check is that $e_p(H(r_1,5))+e_p(H(r_2,5)) <
e_p(H(r_1+r_2,5))$. Thus, we can assume that there is at most one component
equal to $H(r,5)$ and the other components are $K_4$. In fact, replacing
three copies of $K_4$ (contributing $12\cdot 3^p$ to $e_p(G)$ with one copy
of $H(12,5)$ (contributing $11^p+2^{p+1}+9$ to $e_p(G)$) improves $e_p(G)$ so
we can assume that there are at most two components isomorphic to $K_4$.
Since $n \geq 12$ we must have $r \geq 4$. Now, for $r \geq 4$,
$e_p(H(r+4,5)) - e_p(H(r,5)) > 4 \cdot 3^p$ so it is better to replace an
$H(r,5)$ and a $K_4$ with one $H(r+4,5)$. Consequently, $G=H(n,5)$. \endpf

\noindent \textbf{Proof of Theorem \ref{t12}:}\, In the proof we shall
assume, wherever necessary, that $n$ is sufficiently large as a function of
$k$, and that $k \neq 5$. It is trivial to check that the graph $H(n,k)$
defined in the introduction has no $P_k$. We therefore have the lower bound
$t_p(n,P_k) \geq e_p(H(n,k))$. To prove the theorem it suffices to show
that any $P_k$-free graph $G$ with $n$ vertices that is not $H(n,k)$ has
$e_p(G) < e_p(H(n,k))$ for every $p \geq 2$. Assume the contrary, and let
$G=(V,E)$ be a $P_k$-free graph with $n$ vertices that is maximal in the
sense that $e_p(G)=t_p(n,P_k)$ and $G \neq H(n,k)$. We will show how to
derive a contradiction.

According to Corollary \ref{c32}, $|E| \leq n(k-2)/2$. Order the vertices of
$G$ in nonincreasing degree order. That is $V=\{x_1,\ldots,x_n\}$ where 
$d_G(x_i) \geq d_G(x_{i+1})$ for $i=1,\ldots, n-1$. Put $d_i=d_G(x_i)$, and
put $b=\lfloor k/2 \rfloor -1$. Note that
\begin{equation}
\label{e0}
e_p(H(n,k)) = bn^p + o(n^p)
\end{equation}
Put $B=\{x_1,\ldots,x_b\}$ and $A=\{x_{b+1},\ldots,x_n\}$. First observe that
we may assume that $d_1$ is very large. For instance, we may assume that for
all $k \neq 5$, $d_1 > 0.79n$ since otherwise, applying Lemma \ref {l33} to
$G$ with $\alpha=0.79$ and $t=(k-2)/2$ we get for $k \neq 5$, \[ e_p(G) \leq
\frac{k/2-1}{0.79}(0.79)^p n^p +o(n^p) \leq (0.395k-0.79)n^p + o(n^p) <
bn^p+o(n^p) = e_p(H(n,k)). \]

\noindent

\begin{lemma}
\label{l34} If $d_b \leq 0.65n$ then $e_p(G) < e_p(H(n,k))$.
\end{lemma}
\textbf{proof:}\, By (\ref{e0}) it suffices to show that $e_p(G) \leq
cn^p+O(n^{p-1})$ where $c$ is a constant smaller than $b$. Consider the
spanning subgraph of $G$ obtained by deleting all the edges adjacent with the
vertices of $B \setminus \{x_b\}$. Denote this subgraph by $G^{\prime}$. The
maximum degree of $G^{\prime}$ is at most $0.65n$. Let $f_i$ denote the
degree of $x_i$ in $G^{\prime}$ for $i=b,\ldots,n$. By definition,
$e_p(G^{\prime})=f_b^p+\ldots+f_n^p$. Since $f_i \geq d_i-b+1$, and since
$f_b+\ldots+f_n \leq n(k-2) = O(n)$ we have
\begin{equation}
\label{e1}
e_p(G)= d_1^p+\ldots+d_n^p \leq
d_1^p+\ldots+d_{b-1}^p+(f_b+b-1)^p+\ldots+(f_n+b-1)^p =
\end{equation}
$$
d_1^p+\ldots+d_{b-1}^p+e_p(G^{\prime})+o(n^p)
$$
Define $t=e(G^{\prime})/n$. We consider three cases according to the value of
$t$.

\noindent \textbf{Case 1:}\, $t < 0.65$. Since the degree sequence has sum
at most $1.3n$ and no element is larger than $0.65n$ we have
$e_p(G^{\prime}) \leq 2(0.65n)^p$ and using (\ref{e1}) we get: \[ e_p(G) \leq
d_1^p+\ldots+d_{b-1}^p+e_p(G^{\prime})+o(n^p) \leq (b-1+2(0.65)^p)n^p+o(n^p)
< e_p(H(n,k)). \]

\noindent \textbf{Case 2:}\, $1.45 > t \geq 0.65$. According to Lemma \ref
{l33} with $\alpha=0.65$ we know that: \[ e_p(G^{\prime}) \leq
\frac{t}{0.65}(0.65n)^p+ o(n^p) \leq 0.9425n^p+o(n^p). \] Using (\ref{e1}) we
get: \[ e_p(G) \leq d_1^p+\ldots+d_{b-1}^p+e_p(G^{\prime})+o(n^p) \leq
(b-1+0.9425)n^p+o(n^p) < e_p(H(n,k)). \]

\noindent \textbf{Case 3:}\, $t \geq 1.45$. According to Lemma \ref{l33}
with $\alpha=0.65$ we know that:
\begin{equation}
\label{e2}
e_p(G^{\prime}) \leq \frac{t}{0.65}(0.65n)^p+ o(n^p).
\end{equation}
Let $z$ denote the number of edges of $G$ with both endpoints in $B \setminus
\{x_b\}$. Clearly, $z \leq {{b-1} \choose 2} < k^2$. Now, \[ n\frac{k-2}{2}
\geq e(G)=e(G^{\prime})+d_1+\ldots+d_{b-1}-z =tn+d_1+\ldots+d_{b-1}-z \geq
tn+d_1+\ldots+d_{b-1}-k^2 \] It follows that \[ d_1+\ldots+d_{b-1} \leq
n(\frac{k-2}{2}-t)+k^2=n(\frac{k-2}{2}-t)+o(n). \] Since $d_i < n$ the last
inequality immediately gives:
\begin{equation}
\label{e3}
d_1^p+\ldots+d_{b-1}^p \leq (\frac{k-2}{2}-t)n^p+o(n^p).
\end{equation}
Plugging (\ref{e2}) and (\ref{e3}) in (\ref{e1}) yields: \[ e_p(G) \leq
d_1^p+\ldots+d_{b-1}^p+e_p(G^{\prime})+o(n^p) \leq \] \[ (\frac{k-2}{2}-t
+\frac{t}{0.65}(0.65)^p)n^p+o(n^p) \leq (\frac{k-2}{2}-0.5075)n^p+o(n^p) <
e_p(H(n,k)). \] \endpf

In view of Lemma \ref{l34} we may now assume $d_b > 0.65n$, and due to the
remark prior to Lemma \ref{l34} we may also assume that when $k \neq 5$,
$d_1 > 0.79n$. Let $A^{\prime}\subset A$ be the set of vertices that have a
neighbor in $B$. Let $G[A^{\prime}]$ denote the subgraph induced by 
$A^{\prime}$.

\begin{lemma}
\label{l35} If $k$ is even, then $G[A^{\prime}]$ has no edges. If $k$ is
odd, then $G[A^{\prime}]$ contains at most one edge.
\end{lemma}
\textbf{Proof:}\, Assume the contrary. We will derive a contradiction by
showing that $G$ contains a $P_k$. We distinguish three cases

\noindent Consider first the case where $k$ is even. Let $(a_0,a_1)$ be an
edge of $G[A^{\prime}]$. By the definition of $A^{\prime}$, $a_1$ has a
neighbor in $B$. Assume w.l.o.g. $x_1$ is a neighbor of $a_1$. Note that
since $d_i > 0.65n$ for $i=1,\ldots,b$ we have that any two vertices of $B$
have at least $0.3n$ common neighbors in $G$, and hence at least $0.3n-(b-2)
> k$ common neighbors in $A^{\prime}$. Therefore, let $a_i \in A^{\prime}$ be
a common neighbor of $x_{i-1}$ and $x_i$ for $i=2,\ldots,b$ such that
$a_0,a_1,\ldots,a_b$ are all distinct. Let $a_{b+1} \in A^{\prime}$ be a
neighbor of $x_b$ distinct from $a_0,\ldots,a_b$. We have that
$a_0,a_1,x_1,a_2,x_2,a_3,\ldots,a_{b-1},x_{b-1},a_b,x_b,a_{b+1}$ is a $P_k$.

\noindent Next, consider the case where $k$ is odd and there are two edges
in $G[A^{\prime}]$ sharing a common endpoint in $A^{\prime}$. Denote these
two edges by $(a_{-1},a_0)$ and $(a_0,a_1)$. As in the previous case we can
obtain a $P_k$ of the form $a_{-1},a_0,a_1,x_1,a_2,x_2,a_3,
\ldots,a_{b-1},x_{b-1},a_b,x_b,a_{b+1}$.

\noindent Next, consider the case where $k$ is odd and $G[A^{\prime}]$ has
two independent edges, denoted $(a_0,a_1)$ and $(a_{b+1},a_{b+2})$ such that
$a_1$ and $a_{b+1}$ have at least two vertices of $B$ in their neighborhood
union. W.l.o.g. $a_1$ is a neighbor of $x_1$ and $a_{b+1}$ is a neighbor of
$x_b$. As in the previous cases we can obtain a $P_k$ of the form
$a_0,a_1,x_1,a_2,x_2,a_3,\ldots,a_{b-1},x_{b-1},a_b,x_b,a_{b+1},a_{b+2}$.

\noindent The only remaining case is that $k$ is odd and $G[A^{\prime}]$
contains two or more independent edges, and all the endpoints of these
independent edges are connected to a single vertex of $B$, say, $x_1$. In
this case, there may not be a $P_k$ present, but we will show that there is a
$P_k$-free graph $G^{\prime}$ on $n$ vertices with $e_p(G^{\prime}) >
e_p(G)$, contradicting the maximality of $G$. Since we assume $k \geq 7$ we
have $b \geq 2$ so $x_2 \in B$. Let $A^*$ denote the set of non-isolated
vertices in $G[A^{\prime}]$. $|A^*| \geq 4$ and no vertex of $A^*$ is
connected to $x_2$. We may delete the $|A^*|/2$ independent edges of
$G[A^{\prime}]$, and replace them with $|A^*|$ new edges from $x_2$ to each of
the vertices of $A^*$. Clearly, if $G$ is $P_k$-free, so is $G^{\prime}$
(this follows from the fact that $k$ is odd, so $b=\lfloor k/2
\rfloor-1=(k-3)/2$). However, the degree sequence of $G^{\prime}$ majorizes
that of $G$ since the degree of $x_2$ increased, while the other degrees have
not changed. Hence, $e_p(G^{\prime}) > e_p(G)$, a contradiction. \endpf

An immediate corollary of Lemma \ref{l35} is the following:

\begin{coro}
\label{c36} The subgraph of $G$ induced by $B \cup A^{\prime}$ is a spanning
subgraph of $H(b+a^{\prime},k)$ where $|A^{\prime}|=a^{\prime}$. In
particular, if $A^{\prime}=A$ then $G$ is a spanning subgraph of $H(n,k)$.
\endpf
\end{coro}

\noindent Note that by Corollary \ref{c36} we have that if $A^{\prime}=A$
then $e_p(G) < e_p(H(n,k))$ since $G \neq H(n,k)$. This contradicts the
maximality of $e_p(G)$. The only remaining case to consider is when
$A^{\prime}\neq A$. The following lemma shows that this is impossible, due to
the maximality of $G$. This final contradiction completes the proof of
Theorem \ref{t12}.

\begin{lemma}
\label{l37} If $A^{\prime}\neq A$ then there exists a $P_k$-free graph
$G^{\prime}$ with $n$ vertices such that $e_p(G) < e_p(G^{\prime})$.
\end{lemma}

Put $A^{\prime\prime}= A \setminus A^{\prime}$. We claim that each $v \in
A^{\prime\prime}$ has at most one neighbor in $A^{\prime}$. Indeed, if it had
two neighbors, say $a_0,a_1$ then, as in the previous cases, we can obtain a
$P_k$ of the form $a_0,v,a_1,x_1,a_2,x_2,
\ldots,a_{b-1},x_{b-1},a_b,x_b,a_{b+1}$ (in fact, if $k$ is even this is a
$P_{k+1}$). Since $G[A^{\prime\prime}]$ is $P_k$-free it contains at most
$(k/2-1)a^{\prime\prime}$ edges, where $a^{\prime\prime}=|A^{\prime
\prime}|=a-a^{\prime}$. Hence, it contains a vertex $v$ whose degree is at
most $k-2$. Hence $d_G(v)\leq k-1$. Delete all edges adjacent to $v$ in $G$,
and connect $v$ to all edges of $B$. Denote the new graph by $G^{\prime}$.
Note that $G^{\prime}$ is also $P_k$-free. To see this, note that otherwise,
any $P_k$ in $G^{\prime}$ must contain $v$. Let $x_i \in B$ be a neighbor of
$v$ in such a $P_k$. If $v$ is not an endpoint of the $P_k$ it also contains
another neighbor $x_j \in B$ in the path. Since $x_i$ and $x_j$ have many
common neighbors in $A^{\prime}$ (much more than $k$), let $v^{\prime}\in
A^{\prime}$ be such a common neighbor which is not on the $P_k$ (if $v$ is an
endpoint of the $P_k$ it suffices to take $v^{\prime}\in A^{\prime}$ to be
any neighbor of $x_i$ not on the $P_k$). Replacing $v$ with $v^{\prime}$ on
the $P_k$ we obtain a $P_k$ in $G$, contradicting the assumption. We now show
that $e_p(G^{\prime}) > e_p(G)$. Consider the effect of the transformation
from $G$ to $G^{\prime}$ on the degree sequence. The degrees of the vertices
of $B$ increased by one. The degree of $v$ may have decreased by at most
$k-1-b$. The degrees of the neighbors of $v$ in $G$ have decreased by 1.
Since every vertex of $B$ has degree at least $0.65n$, the total increase in
$e_p(G^{\prime})-e_p(G)$ contributed by the vertices of $B$ is at least \[
b((0.65n+1)^p - (0.65n)^p)=bp(0.65n)^{p-1}+o(n^{p-1}). \] Assuming $k \neq
5$, we know $d_1 > 0.79n$. This implies that $ a^{\prime\prime}< 0.21n$.
This fact, together with Lemma \ref{l35} shows that every vertex of
$A^{\prime\prime}$ has degree at most $0.21n$ in $G$. Thus, the total
decrease in $e_p(G^{\prime}) - e_p(G)$ contributed by the vertices of $A$ is
at most \[ (k-1)((0.21n)^p-(0.21n-1)^p)+(k-1)^2-b^2=(k-1)p(0.21n)^{p-1} +
o(n^{p-1}) \] Hence, for $k \neq 5$ \[ e_p(G^{\prime}) - e_p(G) \geq
p(b(3.09)^{p-1}-k+1)(0.21)^{p-1}n^{p-1} + o(n^{p-1}) > 0. \] \endpf

\section{Other acyclic graphs}

A linear forest is a forest whose components are paths. An even linear forest
is a forest whose components are paths with an even number of vertices
(distinct components may have different lengths). The simplest example of an
even linear forest is a \emph{matching}, namely, a graph whose components are
single edges. Let $M_k$ denote the matching with $2k$ vertices. Note that
every even linear forest $F$ with $2k$ vertices is a spanning subgraph of
$P_{2k}$ and contains $M_k$ as a spanning subgraph. Thus, for every $n$ we
have $t_p(n,M_k) \leq t_p(n,F) \leq t_p(n,P_{2k})$. We immediately get the
following proposition:

\begin{prop}
\label{p41} Let $k \geq 2$ be an integer, and let $p \geq 2$ be an integer.
If $F$ is an even linear forest with $2k$ vertices then, for $n$ sufficiently
large, $t_p(n,F)=e_p(H(n,2k))$ where $H(n,2k)$ is the extremal graph
appearing in Theorem \ref{t12}.
\end{prop}
\textbf{Proof:}\, By Theorem \ref{t12} we know that for $n$ sufficiently
large, $t_p(n,P_{2k})=e_p(H(n,2k))$. On the other hand, it is trivial to
check that $H(n,2k)$ does not contain $M_k$ as a subgraph. Hence,
$t_p(n,M_k) \geq e_p(H(n,2k))$. Since $t_p(n,M_k) \leq t_p(n,F) \leq
t_p(n,P_{2k})$ we must have $t_p(n,M_k) = t_p(n,F) = t_p(n,P_{2k})$ for $n$
sufficiently large. \endpf

Another family of trees for which $t_p$ is easy to compute is the family of
stars. Indeed, let $S_k$ denote the star with $k \geq 2$ vertices. Clearly,
if $G$ has no $S_k$ it has $\Delta(G) \leq k-2$. Thus, every $n$-vertex graph
$G$ that is $k-2$-regular must satisfy $t_p(n,S_k)=e_p(G)$. If $n>k-2$ is
even then it is well-known that such $G$ exist for all $k \geq 2$. (in fact,
they can be obtained by an edge-disjoint union of $k-2$ perfect matchings).
So is the case when $n$ is odd and $k$ is even (they can be obtained by an
edge-disjoint union of $(k-2)/2$ Hamilton cycles). If both $n$ and $k$ are
odd then there do not exist $k-2$-regular $n$-vertex graphs, so, clearly, if
$G$ has $n-1$ vertices with degree $k-2$ and one vertex with degree $k-3$,
then $t_p(n,S_k)=e_p(G)$. Such $G$ are well-known to exist for all $n >k-2$.
In fact, they can be obtained by an edge-disjoint union of $(k-3)/2$
Hamilton cycles plus a maximum matching. Note that if $n \leq k-2$, then,
clearly, $t_p(n,S_k)=e_p(K_n)$. To summarize:

\begin{prop}
\label{p42} Let $S_k$ be the star with $k \geq 2$ vertices. Then:

\begin{enumerate}
\item If $n \leq k-2$ then $t_p(n,S_k)=n(n-1)^p$.

\item If $n >k-2$ and $nk$ is even then $t_p(n,S_k)=n(k-2)^p$.

\item If $n > k-2$ and $nk$ is odd then $t_p(n,S_k)=(n-1)(k-2)^p+(k-3)^p$.
\endpf
\end{enumerate}
\end{prop}

A slight modification of $S_k$ is the \emph{near star} $S^*_k$. This graph is
an $S_{k-1}$ to which we add one new neighbor to one of the leaves. So, e.g.,
$S^*_4=P_4$. This slight modification to $S_k$ yields an entirely different
result for $t_p(n,S^*_k)$.

\begin{prop}
\label{p43} If $n > 2k$ then $t_p(n,S^*_k)=e_p(S_n)=(n-1)^p+(n-1)$.
\end{prop}
\textbf{Proof:}\, Let $G$ be a graph without an $S^*_k$. If $G$ has a vertex
of degree $m \geq k-1$ then, trivially, this vertex belongs to a component of
$G$ that is an $S_{m+1}$, since otherwise $G$ would have an $S^*_k$. Hence,
each component of $G$ either has maximum degree at most $k-2$, or else is a
star. Let $s$ denote the number of vertices of $G$ that belong to components
of the first type. Then, $e_p(G) \leq s(k-2)^p+(n-s-1)^p+(n-s-1)$. Clearly,
when $n > 2k$ (in fact, even before that point as $p$ increases), the last
inequality is optimized when $s=0$. Thus, $e_p(G) \leq (n-1)^p+(n-1) $.
Equality is obtained since $S_n$ is $S^*_k$-free. \endpf

A connected bipartite graph is \emph{equipartite} if the two vertex classes
forming the bipartition have equal size. For equipartite trees $T$ that obey
the Erd\H{o}s-S\'os Conjecture we can asymptotically determine $t_p(n,H)$.
Examples of such trees are even paths (however, for these we already have the
sharp result of Theorem \ref{t12}), but there are many others. One example is
the balanced double star $S_{k,k}$, that is obtained by taking two disjoint
copies of the star $S_k$ and joining their roots with an edge. Sidorenko
\cite{Si} has proved that the Tur\'an number of $S_{k,k}$ satisfies
$t(n,S_{k,k}) \leq (k-1)n$ (equality is obtained when $2k-1$ divides $n$).
Namely, the Erd\H{o}s-S\'os Conjecture holds for $S_{k,k}$.

\begin{prop}
\label{p44} If $H$ is an equipartite tree with $2k$ vertices, and $t(n,H)
\leq (k-1)n$ then \[ t_p(n,H) = (k-1)n^p+o(n^p). \]
\end{prop}
\textbf{Proof:}\, We use Lemma \ref{l33} with $\alpha=1$ and $t=k-1$. Indeed,
if $G$ is an $n$-vertex graph that is $H$-free, then $G$ has at most $tn$
edges. Thus, by Lemma \ref{l33}, $e_p(G) \leq (k-1)n^p+o(n^p)$. Consequently,
$t_p(n,H) \leq (k-1)n^p+o(n^p)$. On the other hand, consider the complete
bipartite graph $B_{k-1,n-k+1}$. Since $H$ is equipartite, $B_{k-1,n-k+1}$
does not contain $H$ as a subgraph. Since 
$e_p(B_{k-1,n-k+1})=(k-1)(n-k+1)^p+(n-k+1)(k-1)^p = (k-1)n^p+o(n^p)$ we have
$t_p(n,H) = (k-1)n^p+o(n^p)$. \endpf

\section{Even cycles and complete bipartite graphs}
\textbf{Proof of Theorem \ref{t13}}\, Let $G$ have $n$ vertices and no even
cycle, and assume that $G \neq F_n$. We must show $e_2(G) < e_2(F_n)$. Let 
$d_1 \geq d_2 \geq \ldots \geq d_n$ be the degree sequence of $G$, and let
$x_1,\ldots,x_n$ be the corresponding vertices. Notice first that $d_1+d_2
\leq n+1$. Indeed, otherwise $x_1$ and $x_2$ would have two distinct common
neighbors, and $G$ would contain a $C_4$.

\noindent We first consider the case $d_1 \leq 0.75n$. Straightforward
convexity arguments, plus the fact that $d_1+d_2\leq n+1$ and the fact that
$e(G) < 1.5n$ show that the largest possible value for the sum of squares
is at most the one given by a sequence of the form: \[ 0.75n \; , \;
0.25n(1+o(1)) \; , \; 0.25n(1+o(1)) \; , \; 0.25n(1+o(1) \; , \; d_5 \; , \;
d_6 \; ,\; \ldots\; , \; d_n \] where $d_i$ is bounded by the constant 4 for
$i \geq 5$. Hence $e_2(G) \leq \frac{12}{16}n^2+o(n^2) < e_2(F_n)$ for $n$
sufficiently large (in fact, $n=12$ already suffices).

\noindent Next, we consider the case $d_1 > 0.75n$. Consider any
nonincreasing sequence of $n$ nonnegative integers having the following
properties:

\begin{enumerate}
\item $0.75n < d_1 \leq n-1$.

\item $d_1+d_2 \leq n+1$.

\item There are at least $d_1+1$ nonzero elements in the sequence.

\item The sum of the elements is at most $3(n-1)$.
\end{enumerate}

Putting $d_1=x$, the degree sequence dominates the sequence 
$S=\{x,1,\ldots,1,0,\ldots,0\}$ (there are $x$ ones and $n-x-1$ zeroes
here). Hence, there are at most $3(n-1)-2x$ units to assign to $S$ (subject
to the four properties above) in order to obtain the degree sequence of $G$.
By convexity, the sum of squares is maximized if we assign $n-x$ additional
units to the second, third, etc. elements of $S$, until we run out of units.
Thus, \[ e_2(G) \leq
x^2+(n+1-x)^2\frac{3(n-1)-2x}{n-x}+1^2(x-\frac{3(n-1)-2x}{n-x}). \] Putting
$x=n-k$ (where $1 \leq k < n/4$) the r.h.s. of the last inequality is equal
to $n^2-(k-3)n+3k^2-6$. Thus, \[ e_2(G) \leq n^2-(k-3)n+3k^2-6. \] Note that
when $n > \max\{20 \; , \; 4k\}$ and $k \geq 2$ we have $n^2-(k-3)n+3k^2-6
< n^2+2n-6 \leq e_2(F_n)$. Thus, we have shown that if $d_1 < n-1$ then
$e(G) <e_2(F_n)$. If $d_1=n-1$ then, subject to the above four properties,
the sum of squares is maximized by the unique sequence $n-1,2,2,\ldots,2$.
When $n$ is odd there is only one graph with this degree sequence, namely
$F_n$, and we assume $G \neq F_n$ so $e_2(G) < e_2(F_n)$. When $n$ is even
this is not the degree sequence of any graph (as the sum of the elements is
odd), thus, subject to the above four properties and the requirement that the
sequence be graphic, the sum of squares is maximized by the sequence
$n-1,2,2,\ldots,2,1$. There is only one graph with this degree sequence,
namely $F_n$. Again by our assumption, $G \neq F_n$, so $e_2(G) < e_2(F_n)$.
\endpf

We now turn our attention to complete bipartite graphs. The Tur\'an number
for $K_{k,k}$ is well-understood only for $k=2$. It is known that 
$t(K_{2,2},n) =0.5n^{3/2}(1+o(1))$ (cf. \cite{Bo}). Exact values and extremal
graphs are known only in special cases. Recently, F\"uredi proved in \cite
{Fu} that if a graph has $q^2+q+1$ vertices $q>13$, $m$ edges and no 
$K_{2,2} $ then $m\leq 0.5q(q+1)^2$, and equality holds for graphs obtained
from finite projective planes with polarities. If $k \geq 3$ the asymptotic
behavior of $K_{k,k}$ is not known. The best (and rather simple) bounds are
$t(K_{k,k},n) \leq O(n^{2-1/k})$ and $t(K_{k,k},n) \geq \Omega(n^{2-2/k})$.
Our next proposition shows that $t_k(n,K_{k,k})$ can be asymptotically
determined for every $k \geq 2$. In fact, something slightly stronger can be
proved:

\begin{prop}
\label{p52} Let $2 \leq a \leq k$ where $a,k$ are integers. Then: 
$t_k(n,K_{a,k})=(a-1)n^k(1+o(1))$. Furthermore, if $p \geq k$ then $
t_p(n,K_{2,k})=n^p(1+o(1))$.
\end{prop}
\textbf{Proof:}\, The lower bound is obtained by considering the complete
bipartite graph $K_{a-1,n-a+1}$. It contains no $K_{a,k}$ and obviously has
$e_p(K_{a-1,n-a+1})=(a-1)n^p(1+o(1))$. The upper bound is obtained as
follows. If $G$ has $n$ vertices and no $K_{a,k}$ then \[ \sum_{i=1}^n {{d_i}
\choose k} \leq (a-1){n \choose k} \] since otherwise, by the pigeonhole
principle, there would be at least $a$ vertices whose neighborhood
intersection contains at least $k$ vertices, and hence there would be a
$K_{a,k}$. It follows that $e_k(G) \leq (a-1)n^k(1+o(1))$. In case $a=2$ we
get $e_k(G) \leq n^k(1+o(1))$, so by a trivial convexity argument we get that
for $p \geq k$ $e_p(G) \leq n^p(1+o(1))$. \endpf

\section{Concluding remarks and open problems}

Recently, Pikhurko \cite{Pi} proved the following theorem that was conjectured in the original
version of this paper.
\begin{theo}
\label{t61}[Pikhurko \cite{Pi}]
Let $H$ be a graph with chromatic number $r\geq 3$. Then,
$t_{p}(n,H) = t_p(n,K_r)(1+o(1))$.
\end{theo}

Computing $t_p(n,H)$, or even $t_2(n,H)$ for some specific fixed graphs $H$
seems an interesting open problem. The smallest graph for which we have no
exact answer is $C_4$.

\begin{prob}
\label{p63} Determine $t_2(n,C_4)$. In particular, is it true that for
infinitely many $n$, $t_2(n,C_4)=e_2(F_n)$ where $F_n$ is the friendship
graph.
\end{prob}

Recall that by Proposition \ref{p52}, $t_2(n,C_4)=n^2(1+o(1))$. From Tur\'an
Theory we know that the Tur\'an number of $C_{2k}$ for $k > 2$ (cf.
\cite{Bo}) is smaller than that of $C_4$. This is not the case for $p \geq
2$, since $C_{2k}$ contains $P_{2k}$ and hence $t_p(n,C_{2k}) \geq
(k-1)n^p(1+o(1))$. On the other hand, since $C_{2k}$ is a subgraph of
$K_{k,k}$ we know by Proposition \ref{p52} that
$t_k(n,C_{2k})=(k-1)n^k(1+o(1))$. It is interesting to determine what happens
for $p \neq k$. We conjecture:

\begin{conj}
\label{c64} For $p > 1$, $t_p(n,C_{2k})=(k-1)n^p(1+o(1))$.
\end{conj}

\end{document}